\documentclass[12pt,reqno]{amsart}
\usepackage[margin=1.2in]{geometry}
\usepackage{xcolor}

\usepackage{amsthm, amsfonts, mathrsfs} 
\usepackage{ dsfont } 
\usepackage{amsmath,amssymb} 
\usepackage{enumerate}
\usepackage{graphicx}
\usepackage{subfig}

\usepackage{float}
\usepackage{bbm}
\usepackage{comment}
\usepackage{listings}
\usepackage{hyperref}
\usepackage[sort,compress]{cite}
\usepackage{cite}

\usepackage{color}

\newtheorem{theorem}{Theorem}[section]
\newtheorem{proposition}[theorem]{Proposition}

\newtheorem{definition}[theorem]{Definition}
\newtheorem{remark}[theorem]{Remark}

\newtheorem{conjecture}{Conjecture}

\newtheorem{OP}{Open problem}

\definecolor{bmcolor}{rgb}{0.9, 0.3, 0}

\newcommand{\HH}{\mathcal{H}}

\newcommand{\R}{\mathbb{R}}
\newcommand{\Z}{\mathbb{Z}}
\newcommand{\N}{\mathbb{N}}
\newcommand{\CC}{\mathbb{C}}

\providecommand{\norm}[1]{\lVert#1\rVert}
\newcommand{\G}{\mathcal{G}}
\newcommand{\T}{\mathcal{T}}
\newcommand {\D} {\mathbb D}

\newcommand{\Ww}{\mathcal{W}}
\newcommand{\newmphi}{\kappa_\phi}
\newcommand{\BB}{\mathcal{B}}
\newcommand{\I}{I}
\newcommand{\E}{F}

\newcommand{\supp}{{\rm supp\,}}

\newcommand{\li}{\ell_i}
\newcommand{\mm}{\mathfrak m}
\newcommand{\MM}{\mathfrak M}

\newcounter{lst}

\begin{document}
	
	\date{\today}

\title[Dynamical Sampling]{Dynamical Sampling: A Survey}
\author[A. Aldroubi, C. Cabrelli, I. Krishtal, U. Molter] {Akram Aldroubi, Carlos Cabrelli, Ilya Krishtal and Ursula Molter}

\begin{abstract}
Dynamical sampling refers to a class of problems in which space–time samples are taken from a signal evolving under an underlying dynamical system. The goal is to use these samples to recover relevant information about the system—such as the initial state, the evolution operator, or the sources and sinks driving the dynamics. These problems are tightly connected to frame theory, operator theory, functional analysis, and other foundational areas of mathematics; they also give rise to new theoretical questions and have applications across engineering and the sciences. This survey provides an overview of the theoretical underpinnings of dynamical sampling, summarizes recent results, and outlines directions for future work, including open problems and conjectures.
\end{abstract}

\address{\textrm{(Akram Aldroubi)}
Department of Mathematics,
Vanderbilt University,
Nashville, TN, 37240, USA}
\email{akram.aldroubi@vanderbilt.edu}

\address{\textrm{(Carlos Cabrelli and Ursula Molter)}
	Department of Mathematics,
	UBA and IMAS, UBA-CONICET}
\email{carlos.cabrelli@gmail.com, umolter@dm.uba.ar}

\address{\textrm{(Ilya Krishtal)}
        Department of Mathematical Sciences, Northern Illinois University}
\email{ikrishtal@niu.edu}

\dedicatory{To Ingrid, with deep appreciation of her enormous impact and unparalleled inspiration.}
	
\maketitle


\setcounter{tocdepth}{2}

\section{Introduction and Motivation}

Many natural and engineered processes exhibit dynamics that can be effectively modeled by dynamical systems. Examples include the evolution of temperature fields, the diffusion of pollutants, and the spread of viruses. Understanding and analyzing these evolving phenomena often requires capturing their behavior through measurements taken over time and space. Dynamical sampling, introduced in \cite{ADK13} and inspired by \cite{LV09}, provides a mathematical framework for studying how to efficiently sample and reconstruct time-evolving
signals governed by such models. The term \emph{dynamical sampling} reflects the fact that one samples the solution of a dynamical system in both space and time to extract information about the underlying system. As in most sampling-reconstruction settings, the connection to \emph{frame theory} is central \cite{DS52,AG01, ABK08}.  
Additionally, this area draws heavily on mathematical disciplines such as 
applied harmonic analysis (see, e.g.,~\cite{ACM25, AHKK23, BHKLL25, HLTY25}), 
abstract harmonic analysis (see, e.g.,~\cite{AGHJKR21, KM25, KP25}), 
functional analysis (see, e.g.,~\cite{ADGMM23, Men22, CMS23}), 
and operator theory (see, e.g.,~\cite{ACCMP17, ADK13, ADK15}).
Certain aspects of dynamical sampling are also closely connected to control theory \cite{DMMM2021-1}.

 Dynamical sampling was inspired by problems in applications such as super resolution, sampling and reconstruction in atmospheric emission, source localization of diffusive sources, and distributed sampling of time-varying signals on networks and graphs \cite{LV09,RDCV12,HRLV10, AD18,AD20,MBD15,RBD21,MD17, HLTY25}. The area is closely related to mobile sampling, where a small number of sensors are assumed to move along prescribed trajectories \cite{EG23, GRUV15, GB20, JM22, JMV24, Ku17, RU23, UV13}. 
Beyond its original formulation, dynamical sampling has inspired developments in related applied and theoretical areas. For example, \cite{AMJMP23, AP23} study dynamical dual frames with applications to quantization, \cite{AKT17, AKT18, HM26, LLZ21, Boz24, BK23} consider various questions related to phase retrieval with dynamical frames, \cite{HNT24, AADP13,  ABKMP24, AAP25, LS22} investigate dynamical sampling on structured function spaces, graphs, and other general domains, \cite{Men22, MT23, SCM25, FMT25,FMT25_2, AHP19, Phi17, DLW24, CHP20, SS24} develop  related operator-theoretical techniques for various classes of linear operators in different classes of Hilbert and Banach spaces, \cite{AAK19, AKh17, KM25} perform stability, optimality, and/or perturbation analysis, and \cite{KKT24, KP25} develop techniques for data-driven learning of operators.

Research in dynamical sampling has been shaped to a large degree by the work of Ingrid Daubechies. Indeed, 
dynamical sampling is closely related to Ingrid’s work because it addresses, from a complementary perspective, the same core question that underlies wavelet theory and modern harmonic analysis: how to represent, reconstruct, and stably analyze signals using structured families generated by simple operations. Daubechies’ seminal contributions, such as \cite{Dau92, Dau88, CDF92, DGM86}, showed how stable and efficient representations can be obtained from translates and dilates of a single function, leading to frames and wavelet bases with remarkable localization and robustness properties. Dynamical sampling replaces dilations and translations by the repeated action of an evolution operator (such as a shift, diffusion, or more general linear dynamics), and asks under what conditions time-space samples form a frame or allow stable reconstruction. In this sense, dynamical sampling extends the philosophy of Daubechies’ work -- from static multiscale representations to time-evolving systems -- while relying on the same foundational concepts of frames, stability, redundancy, and operator-based constructions that she helped establish.

There are two related dynamical systems that have often been used as models in dynamical sampling: the discrete-time and continuous-time dynamical systems. 

The discrete-time model is given by
\begin{equation} \label {DS}
\begin{cases}
x(n+1) \;=\; E\,x(n) + w(n), & w(n)\in\Ww,\\[2pt]
y_{n,j} \;=\; \langle x(n), g_j \rangle, & j\in J,
\end{cases}
\end{equation}
where $n$ denotes discrete time, $E$ is a bounded operator on some Hilbert space $\HH$, referred to as the evolution operator, and $\{g_j\}_{j\in J}\subset\HH$ are sampling vectors (sensors); the quantities $y_{n,j}$ are the corresponding space-time measurements, commonly referred to as dynamical samples. 

The continuous-time dynamical system describes the time evolution of a phenomenon with $t\in\R_+$. It is given by
\begin{equation}\label{DFM}
\begin{cases}
\dot{u}(t)=\mathcal A\,u(t)+f(t)+\eta(t),\, t\in\R_+,\\[2pt]
u(0)=u_0,\, u_0\in\HH,\\[2pt]
y_{j}(t_i)=\langle u(\cdot,t_i),\, g_j\rangle,\quad j\in J,\ i\in I,
\end{cases}
\end{equation}
where 
$\HH$ is a Hilbert space of functions on a subset of $\R^d$, such as $L^2([0,1]^d)$. Tthe function $u$ in \eqref{DFM} is $\HH$-valued, and 
$\mathcal A$ is a generator of a $C_0$-semigroup of operators in $B(\HH)$. The space-time samples $y_j(t_i)=\langle u(\cdot,t_i), g_j\rangle$ are obtained using test functions (sensors) $\{g_j\}_{j\in J}\subset\HH$ that are localized in space at locations indexed by $j$.

Dynamical sampling addresses both theoretical and practical challenges in
recovering quantities of interest in dynamical systems from limited
observations. The core problems in this area can be broadly categorized as
\emph{initial state recovery}, \emph{system identification}, and
\emph{source recovery}. These problems have shaped the main research
directions in the field and motivate the organization of this survey.

The fundamental problem of initial state recovery concerns reconstructing
the unknown initial condition from space--time samples of the evolving
system. This line of research began with the pioneering works
\cite{ADK13,ADK15}, where problem \eqref{DS} was studied in the homogeneous
case $\mathcal{W}=0$, with sensing vectors drawn from finite or countable
orthonormal bases and evolution operators given by discrete convolutions.
The complete characterization obtained in \cite{ACMT17} revealed that
initial state recovery is equivalent to the existence of frames generated
by operator orbits. Subsequent work extended these ideas in several directions. One line of
research focuses on general operators and infinite-dimensional settings,
establishing existence, stability, and structural characterizations of
dynamical frames
\cite{CMPP17,CMPP20,AP17,AHP19,ACNP25,ZLL17,MMO23}.

A complementary body of work investigates structured sensing models,
including minimal sensor configurations, graph-based and distributed
settings, and learning-oriented formulations
\cite{AKW15,AADP13,AT14,ACCP21,ACCMP17,AHKL17,
ABKMP24,HLTY25,GH25}.

Source recovery problem addresses the inhomogeneous case $\mathcal{W}\neq 0$ of \eqref{DS}, where the goal is to recover unknown forcing terms driving the dynamics. Early inverse-source formulations for diffusion-type equations appear in \cite{RCLV11,EH02,SW73,LDV11}, while subsequent work developed frame-based and operator-theoretic approaches for discrete-time and continuous-time models
\cite{AD18,MD17,MBD15,RBD21,ADGMM23,ADM24,ACM25,RL25}.
Continuous-time formulations have also been studied using tools from operator semigroup theory and abstract differential equations
\cite{AHKK23,AGK23,AGKMT25}.

A third fundamental problem is system identification, which aims to estimate the evolution operator itself, or its selected spectral and dynamical characteristics, from space-time samples. Representative results utilize generalized Krylov subspaces and other operator- and frame-theoretic  
approaches 
to develop Prony-type methods \cite{AK16,Tan17, KP25}, with further refinements addressing robustness,
partial observability, and numerical aspects
\cite{AHKLLV18,CT22}.

The aim of this paper is to provide a coherent survey of dynamical sampling and its links to frame theory, to present a selection of recent results across continuous and discrete settings (including space-time trade-offs, system identification, and source recovery), and to identify open questions and promising directions for future research.

\subsection{Paper outline and roadmap}
Section~\ref{sec-2} introduces the fundamental problem of initial state recovery and develops the theory of frames of iterations that it motivates. Subsection~\ref{s:ddf_fds} treats dynamical frames in the finite dimensional discrete-time setting: Subsection~\ref{DT} considers frames generated by diagonalizable operators, and Subsection~\ref{GLT} covers the general case. Subsection~\ref{s:ddf_idc} addresses the infinite dimensional discrete-time setting: first the normal operator case in Subsection~\ref{Normal operators}, then the general case in Subsection~\ref{General operators}. Subsections~\ref{Frame-index} and~\ref{Set of generators} study the minimal number of generators needed to obtain a dynamical frame for a given operator. We end with  Subsection~\ref{CTC} which analyzes the continuous-time case.

Section~\ref{SR} investigates recovery of source terms from dynamical samples. Subsection~\ref{DTSR} describes the discrete time source recovery problem, and Subsection~\ref{CTSR} analyzes the continuous time source recovery problem.

Section~\ref{SI} addresses system identification in dynamical sampling, and Section~\ref{outlook} presents open problems and an outlook.

 \begin{remark} We tried to keep the notation close to the original work in which they were introduced. As a consequence, notation may vary across sections.
 \end{remark}

\section{Space-time trade-off and frames of iterations} \label{sec-2}
\subsection{Discrete-time case}
In the discrete-time setting, we assume that an unknown signal evolving in time must be recovered in a stable way from measurements given by sparsely located sensors, which provide local samples at each time step. The idea is to compensate for the limited number of measurements at initial time by exploiting the temporal evolution of the signal, thereby enabling its complete recovery. 

Formally, the framework is as follows: we consider a separable Hilbert space $\HH$, a 
bounded operator $E$ (the evolution operator), and a collection of vectors
$\{g_j\}_{j\in J} \subset \HH$ (the sensors). The goal is to recover every $f \in \HH$ 
in a stable manner from the measurements 
\[
y_{n,j} = \langle E^n f, g_j \rangle, \quad j \in J, \; n = 0,1,...,\ell, \qquad  \text{here } \ell \in \N \cup \{+\infty\}. 
\]
By {\em stable} we mean that if the measurements are not exact but very close, then the resulting $\tilde{f}$ is {\em close} to $f$, i.e. $\|f - \tilde{f}\|$ is small. 

The early dynamical sampling papers \cite{ADK13, ADK15} dealt primarily with the case when $E$ was a convolution operator, i.e.~an operator diagonalizable by the Fourier transform. The dynamical samples were typically taken on a spatial grid and could be written in the form
\begin{equation}
    \label{e:sds}
    y_n = S_\Omega E^n f,\ 
\end{equation}
where $S_\Omega$ is an orthogonal projection onto the subspace $\HH_\Omega = \overline{\rm span}\, \{e_j: j\in\Omega\}$ for a fixed orthonormal basis $\{e_j\}_j$. Based on the classical results of Papoulis \cite{P77}, the papers \cite{ADK13, ADK15} characterized structured sets $\Omega$ that yield a positive solution of the time-space trade-off dynamical sampling problem and obtained numerical estimates for stability. 

The complete solution of the problem, obtained in \cite{ACMT17}, is based on the following simple but crucial observation:  
\begin{equation}\label{e:trick}
y_{n,j} =\langle E^n f, g_j \rangle = \langle f, (E^*)^n g_j \rangle.    
\end{equation}
Hence, the exact condition for the stable recovery of $f$ is that the system 
$\{(E^*)^n g_j : j \in J,\, n = 0,1,...,\ell\}$ forms a frame for $\HH$. That is, there exist positive constants $\mm$ and $\MM$, such that 
$$ \forall f \in \HH, \quad \mathfrak m \|f\|^2 \leq \sum_{j \in J,\, 0\leq n \leq \ell} \left |\langle f, (E^*)^n g_j \rangle \right |^2 \leq \mathfrak M \|f\|^2. $$

Thus, we get the following equivalent problem: 
given a bounded operator $T$ acting on a separable Hilbert space $\HH$, 
find conditions on $T$ 
and a collection of vectors 
$\{g_j\}_{j \in J} \subset \HH$ such that the family 
\begin{equation} \label{frame}
\{T^n g_j : j \in J,\, n = 0,1,...,\ell\}
\end{equation}
forms a frame for $\HH$. Such frames are called {\it frames of iterations} or {\it dynamical frames}.

We note that observation \eqref{e:trick} plays a significant role not only in the study of the time-space trade-off problem but also in many other dynamical sampling problems.

\subsection{Discrete-time dynamical frames. Finite dimensional case}\label{s:ddf_fds}
For a $d$-dimensional complex Hilbert space $\HH \simeq \CC^d$, any set of generators is a frame. So the question becomes: Under what conditions on the operator $T$ acting on $\CC^d$ and vectors $\{b_j\}_{j\in J} \subset \CC^d$, is the set $\{T^nb_j\}_{j \in J,\, n = 0,1,...,\ell}$  a set of generators of $\CC^d$.

In this case, there is a complete characterization of how many vectors need to be iterated and how many iterations are necessary to yield a frame.

With a slight abuse of notation, we identify the operator $T$ with a $d\times d$-matrix in $\CC^{d\times d}$. We can reduce the problem as follows: Let $B \in \CC^{d\times d}$ be any invertible matrix, and let $Q$ be the matrix $Q=BTB^{-1}$, i.e.\ $T = B^{-1}QB$. Let also $b_i$ denote the $i$-th column of $B$.
Since a frame is transformed into a frame by invertible linear operators, condition (\ref{frame})  is equivalent to $\{Q^jb_i: i \in \Omega, \; j=0,\dots, \li\}$ being a frame of $\CC^d$.

Since any matrix $T$ admits a {\em Jordan canonical form}, we need to study only matrices in Jordan form. Since the treatment is different, we separate two cases.

\subsubsection{Diagonalizable Transformations} \label {DT}
Let $T\in \CC^{d\times d}$ be a diagonalizable matrix, with distinct eigenvalues $\{\lambda_1,\dots,\lambda_n\}$, i.e. $T=B^{-1}DB$ with $D$  a diagonal matrix of the form
\begin {equation}\label {Adiag}
D=
 \begin{pmatrix}
 
  \lambda_1 I_1 & 0 & \cdots & 0 \\
  0 & \lambda_2I_2 & \cdots & 0 \\
  \vdots  & \vdots  & \ddots & \vdots  \\
  0 & 0 & \cdots & \lambda_n I_n
 \end{pmatrix} \text{ with } I_k \text{ an } h_k\times h_k \text{ identity matrix, } k = 1, \dots n. 
 \end {equation}
Recall that the {\em $T$-annihilator $q^T_b$ of a vector $b$} is the monic polynomial of smallest degree, such that $q_b^T(T)b \equiv 0$, and let  $P_j$ denote 
the orthogonal projection in $\CC^d$ onto the eigenspace of $D$ associated to the eigenvalue $\lambda_j$. Then we have:
\begin {theorem}[\cite{ACMT17}]\label {TAD0}
Let $\{b_i: i \in J\}$  be vectors in $\CC^d$. Let $D$ be a diagonal matrix and $r_i$ the degree of the $D$-annihilator of $b_i$. Set $\li=r_i-1$. Then  $\{D^{j}b_i: \; i\in J,  \, j=0, \dots, \li\}$ is a frame of $\CC^d$ if and only if $\{P_j(b_i):i \in \Omega\}$ form a frame of $P_j(\CC^d)$, $j=1,\dots ,n$.
\end{theorem}
This theorem provides a complete characterization of all  subsets of vectors $\{ h_k \}$ for which a frame generated by iterations exists, as well as the number of iterations required. Specifically, for any set $\{ h_k \}$ satisfying the conditions of the theorem above, exactly 
$\ell = \sup_{k = 1, \dots, n} h_k$ iterations are needed to generate a frame using the set $\{ h_k \}$.

An interesting special case of Theorem \ref{TAD0} was already observed in \cite{ADK13}: it was noted that if the evolution operator is a convolution operator with distinct eigenvalues, then every vector can be recovered from dynamical samples at a single node if and only if none of the components of that node, in the basis of eigenvectors, is zero.

\subsubsection {General linear transformations}
\label {GLT}

For a general matrix, in order to state the results in this case, we need to introduce some notation related to the 
general Jordan form of a matrix with complex entries. A matrix $J \in \CC^{d\times d}$ with distinct eigenvalues $\lambda_s, s = 1, \dots, n$  is in Jordan form if 
\begin {equation}\label {AJord}
J=
 \begin{pmatrix}
 
\lambda_1I_1+ N_1 & 0 & \cdots & 0 \\
  0 & \lambda_2I_2+ N_2 & \cdots & 0 \\
  \vdots  & \vdots  & \ddots & \vdots  \\
  0 & 0 & \cdots & \lambda_nI_n+ N_n
 \end{pmatrix}, \ 
 N_{s} =
 \begin{pmatrix}
   N_{s_{1}} & 0 & \cdots & 0 \\
  0 & N_{s_{2}} & \cdots & 0 \\
  \vdots  & \vdots  & \ddots & \vdots  \\
  0 & 0 & \cdots & N_{s_{\gamma_s}},
 \end{pmatrix} 
 \end {equation}
 where, as before $I_{s}$ is an $h_s\times h_s$ identity matrix, $h_1 + \dots + h_n = d$, and each $N_{s_i}$ is a $t^{s}_i \times t^{s}_i$ cyclic nilpotent matrix,  
\begin {equation}\label {cyc}
N_{s_i} = 
\begin{pmatrix}
0 & 0 & 0 & \cdots & 0 & 0\\
1 & 0 & 0 & \cdots & 0 & 0\\
0 & 1 & 0 & \cdots & 0 & 0\\
\vdots  & \vdots & \ddots  & \vdots & \vdots  & \vdots\\
0 & 0 & 0 & \cdots & 1 & 0
\end{pmatrix},
\text{ with } t^{s}_1\geq t^{s}_2 \geq \dots,\text{ and }t^{s}_1 + t^{s}_2 + \dots + t^s_{\gamma_s} = h_s. 
\end {equation}

Let $k_j^s$ denote the index corresponding to the first row of the block $N_{s_j}$ from the matrix $J$, and let $e_{k_j^s}$ be the corresponding element of the standard basis of $\CC^d$. (That is a cyclic vector associated to that block). 
We also define $W_s := \text{span} \{e_{k^s_j}: j=1,\dots, \gamma_s\}$, for $s=1, \dots,n$, and $P_s$ will again denote the orthogonal projection onto $W_s$. Using the notation and definitions above, we can state the following theorem:
\begin {theorem}[\cite{ACMT17}, Theorem 2.6]
\label {tadcul}
Let $T$ be a matrix, such that $T = B^{-1}JB$, where $J \in \CC^{d\times d}$ is the Jordan matrix for $T$. Let $\{b_i: i \in \Omega\}$ be a subset of the column vectors of $B$, $r_i$ be the degree of  the $J$-annihilator of the vector $b_i$, and  let $\li=r_i-1$.\\
Then the following propositions are equivalent.
\begin{enumerate}
\item[i)]
The set of vectors $\{T^jb_i: \; i \in \Omega, j=0,\dots, \li\}$ is a frame for $\CC^d$.
\item[ii)]
For every $s = 1, \dots, n$, \;$\{P_s(b_i), i \in \Omega\}$ form a frame of $W_s.$
\end{enumerate}
\end{theorem}
In other words, $\{T^jb_i: \; i \in \Omega, j=0,\dots, \li\}$ is a frame for $\CC^d$ if and only if the {\em Jordan}-vectors of $T$  (i.e.\ the columns of the matrix $B$ that reduces $T^*$ to its Jordan form) corresponding to $\Omega$ satisfy that their projections onto the spaces $W_s$ form a frame. 

This theorem again provides a complete characterization of all  subsets of vectors $\{ h_k \}$ for which a frame generated by iterations exists, as well as the number of iterations required. In this case, we need at least 
$\ell = \sup_{s=1, \dots, n} \{ s_{\gamma_s}\}$ vectors  $\{h_k, : k = 1, \dots,\ell\}$ that satisfy the conditions of the preceding theorem.

\subsection{Discrete-time dynamical frames. Infinite dimensional case}\label{s:ddf_idc}

Let \(\HH\) be a complex separable Hilbert space.
In this section, we consider \emph{infinite} iterations of an operator, i.e., systems of the form
\(\{T^{n}g\}_{g\in\G,\;n\ge 0}\).
In this regime, we can characterize precisely which operators admit a frame of iterations.
Perhaps surprisingly, many operators do \emph{not} admit any frame of iterations.

Before proceeding, we record a useful sufficient condition for the existence of \emph{tight} dynamical frames.

\begin{theorem}[\cite{AP17}]\label{thm:tight-frame-AP17}
If \(T\) is a contraction (\(\|T\|\le 1\)) and \((T^*)^{n}f\to 0\) as \(n\to\infty\) for every \(f\in\HH\) (strong stability of \(T^*\)),
then there exists a set \(\G\subseteq\HH\) such that \(\{T^{n}g\}_{g\in\G,\;n\ge 0}\) is a \emph{tight} frame for \(\HH\).
\end{theorem}

This motivates the following definitions.

\begin{definition}[Contraction and strong stability]
Let $\mathfrak{U}$ be a bounded operator on a Hilbert space $\HH$.
\begin{itemize}
  \item $\mathfrak{U}$ is a \emph{contraction} if $\|\mathfrak{U}\|\le 1$.
  \item $\mathfrak{U}$ is \emph{strongly stable} if $\mathfrak{U}^{\,n}f \to 0$ for every $f\in\HH$ (i.e., $\mathfrak{U}^{\,n}\!\xrightarrow{\;\text{sot}\;}0$ as $n\to\infty$).
\end{itemize}
\end{definition}

\begin{definition}[Similarity]
Given Hilbert spaces $\HH_1,\HH_2$, operators $T_1\in\BB(\HH_1)$ and $T_2\in\BB(\HH_2)$ are \emph{similar} if there exists a bounded invertible operator
$P:\HH_1\to\HH_2$ such that
\[
P^{-1}T_2P=T_1.
\]
\end{definition}

We are now in position to state the necessary and sufficient conditions for the existence of dynamical frames.

\begin{theorem}[\cite{CMPP20,ACNP25}]\label{thm:dynamical-frame}
Let $\HH$ be a complex separable Hilbert space and $T\in\BB(\HH)$.
Then the following are equivalent:
\begin{enumerate}
  \item There exists a collection $\{g_j\}_{j\in J}\subset\HH$ such that
  $\{T^n g_j : j\in J,\; n\ge 0\}$ is a (dynamical) frame for $\HH$.
  \item $T$ is similar to a contraction $C$ whose adjoint $C^*$ is strongly stable.
\end{enumerate}
In particular, (2) means: there exists a bounded invertible $P$ with $P^{-1}TP=C$, $\|C\|\le 1$, and $(C^*)^{\,n}f\to 0$ for every $f\in\HH$. 
\end{theorem}

We may now ask which operators can generate a {\em Parseval }frame.
In this situation, the operator must be a contraction. Thus, we obtain:

\begin{theorem}[\cite{ACNP25}]\label{thm:dynamical-Parseval-frame}
Let $\HH$ be a complex separable Hilbert space and $T\in\BB(\HH)$.
Then the following are equivalent:
\begin{enumerate}
  \item There exists a collection $\{g_j\}_{j\in J}\subset\HH$ such that
  $\{T^n g_j : j\in J,\; n\ge 0\}$ is a  Parseval frame for $\HH$.
  \item $T$ is a contraction  whose adjoint is strongly stable.
\end{enumerate}
\end{theorem}
Thus, only contractions can produce Parseval frames of iterations!!
Parseval frames are important since its canonical dual is the same frame, what is crucial in applications.

It is known that given a frame $\{g_j\}_{j\in J}$ of $\HH$, then  $\{S^{-1/2}(g_j)\}_{j\in J}$, is a Parseval frame of $\HH$ where $S$ is the frame operator of $\{g_j\}_{j\in J}$.

To see this we use that a frame is Parseval if and only if its frame operator is the identity.
Thus, if $R$ is the frame operator of $\{S^{-1/2}(g_j)\}_{j\in J}$, then 
$$
Rf=\sum_{j\in J} \langle f, S^{-1/2}(g_j) \rangle S^{-1/2}(g_j)= S^{-1/2} (
\sum_{j\in J} \langle S^{-1/2}f, g_j \rangle g_j) = S^{-1/2}S S^{-1/2}f = f.
$$
Applying this to a frame of iterations $\{T^ng_j\}_{n,j}$ we get that $\{S^{-1/2}T^ng_j\}_{n,j}$ is a Parseval frame of iterations since $S^{-1/2}T^ng_j=(S^{-1/2}Tg_jS^{1/2})^n(S^{-1/2}g_j)$.

It is interesting to mention here that the canonical dual of a dynamical frame is also a dynamical frame \cite{AKh17}. This follows immediately  from the fact that  the canonical dual of a frame is the image of the frame by $S^{-1}$ where $S$ is the frame operator, and the following inequalities:
$$
S^{-1}T^ng_j = (S^{-1}T^nS)\, S^{-1}g_j = (S^{-1}TS)^n S^{-1}g_j.
$$

\subsubsection{Normal operators}
\label{Normal operators}
The case of a normal operator is particularly rich, establishing deep connections with several areas of mathematics, including complex analysis, Hardy spaces, and spectral theory. It was among the first cases to be studied, and to this day it still presents fascinating open problems. Moreover, it remains one of the few situations where explicit examples of dynamical frames are available.

The first surprising result is that there does not exist a Riesz basis of iterations for a normal operator:

\begin{theorem}[\cite{ACMT17,ACCMP17,AP17}]
If  $T\in\BB(\HH)$ is  a normal operator. Then for any collection of vectors 
$\G\subset\HH$ the family of iterates $\{T^ng\}_{g\in \G, n\geq 0}$ is not a basis of $\HH.$
\end{theorem}

This follows as a consequence of the Müntz–Szász theorem (see \cite{ACMT17}).  
Therefore, our next step is to investigate the existence of frames.  
Recall that strongly stable operators play a key role in this context.  
In particular, since unitary operators are never strongly stable, Theorem~\ref{thm:dynamical-frame} implies that they cannot generate frames of iterations.  
We then begin by describing the simple case of diagonal operators acting on $\ell^2(\N)$ in the case where only a single vector is iterated.  

\begin{theorem}\cite{ACMT17} \label {OnePointFrame}
Let  $D=\sum_j\lambda_jP_j$ be such that $P_j$ have rank $1$ for all $j\in \N$, and let $b := \{b(k)\}_{k \in \N} \in \ell^2(\N)$. 
Then $\{D^lb: l=0,1,\dots\}$ is a frame if and only if
\begin{enumerate}
\item[i)] $|\lambda_k| < 1$  for all $k.$  
\item[ii)]  $|\lambda_k| \to 1$.
\item[iii)] $\{\lambda_k\}$ satisfies Carleson's condition 
\begin{equation}
\label{carleson-cond}
\inf_{n} \prod_{k\neq n} \frac{|\lambda_n-\lambda_k|}{|1-\bar{\lambda}_n\lambda_k|}\geq \delta.
\end{equation}
for some $\delta>0$.
\item[iv)] $b(k)=m_k\sqrt{1-|\lambda_k|^2}$ for some sequence $\{m_k\}$ satisfying $0<C_1\le |m_k| \le C_2< \infty$.

\end{enumerate}
\end{theorem}
The condition that the eigenvalues have modulus less than one makes the operator contractive and strongly stable.
Items ii) and iii) are more subtle and connect this problem with interpolation sequences in the Hardy space in the complex open unit disk $\mathbb D$ \cite{ACMT17}.
These frames are typically called Carleson frames (sometimes also Alcamota frames).
They have remarkable features that have not been identified for any other frames in the literature, such as excessive redundancy. In particular, it was observed in \cite{CHPS24} that the subfamily obtained by selecting each $N^{\rm th}$ element from a Carleson frame remains a frame under a mild additional condition on the spectrum of $D$ regardless of the choice of $N\in \N.$
 Furthermore, the frame property is preserved upon removal of an arbitrarily finite number of elements. In \cite{KM26}, the result was extended further to establish that systems of the form $\{D^{Nk+j_k}b\}_{k\in\N}$ are frames for $\ell^2(\N)$ for any choice of $N\in \N$ and $j_k \in [0, N)$ assuming that $D$ is positive definite. The latter result can be proved in two different ways: it relies either on an application of the M\"untz–Sz\'asz theorem or on a deeper result of Rubel \cite{R56} about the density of zeros of entire functions. The full extent of excessive redundancy of Carleson frames remains unknown. We offer the following conjecture.

 \begin{conjecture}\label{conj:carl}
     Let $\{D^lb: l=0,1,\dots\}$ be a Carleson frame with $\sigma(D) \subset (0,1)$. Assume also that a set $\Lambda \subset (0,\infty)$  satisfies the M\"untz–Sz\'asz condition
     \[
     \sum_{\lambda \in \Lambda} \frac1\lambda = \infty.
     \]
     Then the system $\{D^\lambda b: \lambda\in\Lambda\}$ is a frame for $\ell^2(\N)$, as long as it is a Bessel sequence.
 \end{conjecture}

In a series of papers \cite{ACCMP17,AP17,CMPP20,CMS21}, the authors considered frames of iterations generated by finitely many vectors and general normal operators. These works provide necessary conditions on vectors to generate dynamical frames. The main tools come from spectral theory with multiplicity and from function theory in the Hardy space $H^2(\mathbb{D})$ 
-- especially Blaschke products and interpolation sequences -- revealing a rich network of connections among these areas. 

\subsubsection{General operators}
\label {General operators}
In \cite{CHP20}, the authors observed that if a bounded operator \(T\) on a Hilbert space \(\HH\) admits a vector whose orbit 
\(\{T^n g\}_{n \ge 0}\) forms a frame, then the kernel of the synthesis operator associated with this frame 
can be realized as an invariant subspace $M$ for the shift operator acting on the Hardy space \(H^2(\D)\), 
where \(\D\) denotes the open unit disk in the complex plane. 

Consequently, if we denote by \(N\) the orthogonal complement of this kernel in \(H^2(\D)\), 
the restriction of the synthesis operator to \(N\) defines an isomorphism from \(N\) onto \(\HH\) 
that intertwines \(T\) with \(A_N\), the compression of the shift to \(N\). 

This correspondence allows one to study our problem for the operator $A_N$ acting in $N$,  a space of richer analytic structure, 
and then transfer the obtained results back to the original space \(\HH\) and the operator $T$. 
Subspaces of the Hardy space whose orthogonal complement is shift-invariant are known as \emph{model spaces}, 
and they play a central role in the model theory of contractions (see \cite{NagyFoias2010}).
The result in \cite{CHP20} has been generalized to an arbitrary number of generators in
\cite{CMS23}.

The strategy in these works is to construct functional models 
that intertwine the operator under study with simpler operators 
acting on suitable function spaces. In some cases, the 
intertwining map is merely an isomorphism, as in the examples 
mentioned above; however, the framework becomes significantly 
more powerful when this intertwining is realized through a 
\emph{unitary} operator. The results that follow, as well as  the existence results in
Theorem \ref{thm:dynamical-frame} and Theorem \ref{thm:dynamical-Parseval-frame} mentioned before, make extensive 
use of the model introduced in \cite{Rota1960} and later refined 
in \cite{deBrangesRovnyak1964}.

\subsubsection{Frame-index of an operator}
\label{Frame-index}
Once the existence results have been established, that is, once we know precisely which operators admit frames of iterations, an essential issue becomes characterizing the generating collections. A question of particular interest is to determine the smallest number of generators capable of producing a frame of iterations, as this number carries both theoretical significance and practical implications.

In many settings arising from dynamical sampling, the main goal is to recover the initial state of an evolving process from a finite set of observations gathered by sensors fixed at certain spatial points and recording the system at successive time steps. Within this interpretation, sensor locations play the role of generators. From an applied perspective, identifying the smallest collection of generators is crucial, as it corresponds to the most economical configuration of sensors, which still guarantees stable reconstruction of the evolving signal. Reducing the number of generators not only decreases costs and addresses constraints such as accessibility or hardware limitations, but also sheds light on intrinsic structural features of the operator itself. This balance between efficiency and structure provides strong motivation for the study of  minimal generating systems.

The \emph{minimum cardinality} of a set of generators for a given operator \(T\) 
is called the \emph{frame index} of \(T\), and it is denoted by \(\gamma(T)\). 
When the operator \(T\) does not admit any frame of iterations, we set 
\(\gamma(T) = 0\); the index is said to be \emph{infinite} when \(T\) admits 
frames of iterations, but none can be generated by a finite number of vectors. 
The \emph{Parseval index} of an operator, denoted by \(\gamma_p(T)\), is defined 
analogously, replacing ``frame'' by ``Parseval frame'' in the above definition.

The first attempt to estimate the frame index of a bounded operator \(T\) appears in 
\cite[Proposition 2.13]{CMS23}. In that work, the authors studied the problem 
within a model space associated with the compression of the shift (the model 
operator). There is no loss of generality in considering this particular operator, 
since every operator that admits a frame of iterations is known to be similar to a model operator, and the index is 
invariant under similarity. Using deep results from complex analysis and operator 
theory -- such as the Corona Theorem and the Gelfand transform -- they established 
a lower bound for the index, which they conjectured to be sharp.

More recently, in \cite{ACNP25}, the authors  obtain the following results using functional models.
For the Parseval index the following hold.

\begin{theorem}[\cite{ACNP25}]\label{nikc2}
Let $T \in\mathcal{B}(\HH)$. If $T$ has a Parseval frame of iterations, then
\begin{equation*}
\gamma_p(T)=\dim\overline{(I_\HH-TT^*)(\HH)}.
\end{equation*}
Moreover, the index $\gamma_p(T)$ is attained with linearly independent generators.
\end{theorem}

On the other side, for the general index they obtain the following proposition.
\begin{proposition}[\cite{ACNP25}]\label{general index}
Given an operator $T\in\mathcal{B}(\HH)$ we have: 
\begin{align*}
\gamma(T) &=\min\;\{\gamma_p(Q) : 
Q\text{ a strongly stable contraction  similar to } T\}.
\end{align*}
Here we take $\min\{\emptyset\}=0$.
\end{proposition}

When \(T\) is a contraction that is {\em jointly strongly stable} -- that is, both 
\(T\) and \(T^*\) are strongly stable -- then each of these operators admits frames 
of iterations. This case was studied in detail in \cite{ACNP25}, where conditions 
were established to ensure that the frames generated by them are similar.
\subsubsection{Set of generators}
\label {Set of generators}

In \cite{CHP20} it was shown that if a vector $g \in \HH$ is a 
frame generator for an operator $T \in \mathcal{B}(\HH)$, 
then every other generator is of the form $Qh$, 
where $Q$ is a bounded invertible operator that commutes 
with~$T$. Moreover, in \cite{CMS23}, a characterization of all 
finite sets of generators was obtained for a model operator. 
This result can then be transferred to any bounded operator 
admitting a finite set of generators.

Additional contributions to frames generated by operator orbits include commuting and multi-operator constructions 
\cite{ACCP21, ACCP23, ACNP25, BC25}, finite-dimensional and cyclic frame constructions \cite{CRS26, SR25}, multivariate, group, and tensor-product extensions \cite{BHKLL25, HX22, Men22, ZLL17_1, ZL19}, and operator-theoretic characterizations of Bessel and closed-range orbits \cite{Men22}.

\subsection{Continuous-time case: formulation and semi-continuous frames} \label {CTC}
In the continuous-time setting, an initial state \(f\in\HH\) evolves under a bounded operator \(A\in\mathcal{B}(\HH)\) via
\begin{equation} \label{CTEO}
f_t = A^{t}f,\qquad t\ge 0.    
\end{equation}

We observe this evolution through a finite or countable sensing set \(\G\subset\HH\) at times \(t\) in a set \(\I\subset[0,\infty)\), collecting the samples
\[
\bigl\{\langle A^{t}f,\,g\rangle:\ g\in\G,\ t\in \I\bigr\}.
\]

As in the discrete case, the core questions are: for which triples \((A,\G,\I)\) is every \(f\in\HH\) recoverable in a stable way? Equivalently, recovery asks when the family \(\{(A^*)^{t}g:\ g\in\G,\ t\in\I\}\) is complete in \(\HH\); stable recovery corresponds to frame inequalities for these adjoint iterates. For convenience, we replace \(A\) with \(A^*\) and drop the asterisk in what follows; thus we ask when \(\{A^{t}g:\ g\in\G,\ t\in\I\}\)  forms a semi-continuous frame?

When time is an interval \(I=[0,L]\), it is convenient to view measurements over the product set \(\G\times[0,L]\) with the product measure “(counting on \(\G\)) \(\times\) (Lebesgue on \([0,L])\)”. In this case the collection
\[
\{A^{t}g:\ g\in\G,\ 0\le t\le L\}
\]
is a \emph{semi-continuous frame} for \(\HH\) if there exist constants \(0<\mathfrak m\le \mathfrak M<\infty\) such that
\begin{equation}\label{SCF}
  \mathfrak m\|f\|^{2}\ \le\ \sum_{g\in\G}\int_{0}^{L}\bigl|\langle f,\,A^{t}g\rangle\bigr|^{2}\,dt\ \le\ \mathfrak M\|f\|^{2}
  \qquad\text{for all }f\in\HH.
\end{equation}
Informally, \eqref{SCF} captures the continuous-time trade-off: integrating information over time compensates for using fewer spatial sensors. The discrete-time formulation is recovered by replacing the integral with a sum over a finite time mesh \(\T\subset[0,L]\); below we indicate when such discretizations preserve stability.  

\medskip

\begin{theorem}\cite{AHP19} \label{ScToDscr}
Let \(A\in\mathcal{B}(\HH)\) be normal and let \(\G\) be a Bessel family in \(\HH\). If \(\{A^t g \}_{g\in\G,\ t\in[0,L]}\) is a semi-continuous frame for \(\HH\), then there exists \(\delta>0\) such that for any finite set
\(\T=\{t_i\}_{i=1}^{n+1}\subset[0,L]\) with \(0=t_1<\dots<t_{n+1}=L\) and mesh
\(\max_{1\le i\le n}(t_{i+1}-t_i)<\delta\),
the system \(\{A^{t}g\}_{g\in \G,\ t\in \T}\) is a frame for \(\HH\).

If, in addition, \(A\) is invertible, then \(\{A^t g\}_{g\in\G,\ t\in[0,L]}\) is a semi-continuous frame for \(\HH\) \emph{if and only if} there exists a finite set \(\T=\{t_i\}_{i=1}^{n+1}\subset[0,L]\) with \(0=t_1<\dots<t_{n+1}=L\) such that \(\{A^{t}g\}_{g\in \G,\ t\in \T}\) is a frame for \(\HH\).
\end{theorem}

The invertibility of \(A\) is necessary for the equivalence in the second statement of Theorem~\ref{ScToDscr}.

The theorem above links semi-continuous frames generated by the iterates \(\{A^{t}g\}\) to discrete-time counterparts, including frames generated by nonuniform powers of \(A\).
A general framework connecting continuous frames and their discretizations can be found in \cite{FS19}.

In infinite-dimensional Hilbert spaces, Theorem~\ref{ScToDscr} implies that \(|\G|=\infty\): since \(\mathcal{T}\) is chosen to be finite and any frame \(\{A^{t}g\}_{g\in\G,\ t\in\mathcal{T}}\) must have infinitely many vectors,  \(|\G|\) cannot be finite.

Theorem~\ref{ScToDscr} admits several generalizations; see \cite{DMM21}. In particular, one may take the time interval to be \([0,\infty)\) or consider strongly continuous semigroups \(\{e^{tA}\}_{t\ge 0}\) with \(A\) not necessarily normal. For \emph{exponentially stable} semigroups, the same discretization conclusion holds on \([0,\infty)\). In contrast, for semigroups that are \emph{not} exponentially stable, it may happen that \(|\G|\) is finite, but then the time set must be infinite (i.e., discretization requires infinitely many time samples). Note, however, that Theorem~\ref{ScToDscr} is not a special case of the semigroup setting: the two notions \(A^{t}\) (defined via the spectral calculus) and \(e^{tA}\) do not coincide in general.

It is noteworthy that if \(\{A^t g\}_{g\in\G,\, t\in[0,L_0]}\) is a semi-continuous frame for some \(L_0>0\), then it is a semi–continuous frame for every finite \(L>0\). The following theorem makes this precise.

\begin{theorem} \cite{AHP19}\label{SCFrSA}
Let \(A\in\mathcal{B}(\HH)\) be invertible and self-adjoint, and let \(\G\) be a countable subset of \(\HH\).
Then \(\{A^t g\}_{g\in\G,\, t\in[0,1]}\) is a semi-continuous frame for \(\HH\) \emph{if and only if}
\(\{A^t g\}_{g\in\G,\, t\in[0,L]}\) is a semi-continuous frame for \(\HH\) for every finite \(L>0\).
\end{theorem}

\begin{conjecture}\label{conj:reduct}
The conclusion of Theorem~\ref{SCFrSA} remains true when \(A\) is normal and reductive.
\end{conjecture} 
\begin{OP} \label {OP:L-infty}
 What changes in Theorems~\ref{ScToDscr} and~\ref{SCFrSA} if the assumption \(L<\infty\) is replaced by \(L=\infty\)?
 \end{OP}

\subsubsection{Convolution-kernel (semi-group) case}
A central applied setting is when \(\HH\) is the Paley–Wiener space of band-limited functions
\[
PW_{c} \;:=\; \bigl\{ f\in L^{2}(\mathbb{R}) : \operatorname{supp}\widehat{f}\subseteq[-c,c]\bigr\},
\]
and the evolution acts by convolution with a family of kernels \(\{\phi_{t}\}_{t>0}\):
\[
f_t(x) \;=\; (\phi_t * f)(x),\qquad t>0,
\]
with the semi-group properties \(\phi_{t+s}=\phi_t * \phi_s\) for \(t,s>0\) and \(\phi_t*f \to f\) in \(L^2(\mathbb{R})\) as \(t\to 0^+\).
In addition, for physical systems, the set of kernels is assumed to satisfy: 
\begin{equation}\label{DefOp}
\Phi_c:=\bigl\{\phi\in L^1(\mathbb{R})\;\big|\;\exists\,\newmphi>0:\ (\forall\,\xi\in[-c,c])\,(\newmphi\le\widehat{\phi}(\xi)\le 1)\, \ \widehat{\phi}(0)=1\bigr\}.
\end{equation}

A prototypical example is the (one-dimensional) fractional diffusion equation
\[
\begin{cases}
\partial_t u(x,t) \,=\, (-\partial_x^2)^{\alpha/2} u(x,t), & x\in\mathbb{R},\ t>0,\\[2pt]
u(x,0)=f(x),
\end{cases}
\qquad 0<\alpha\le 2,
\]
whose solution is \(u(\cdot,t)=\phi_t * f\) with
\[
\widehat{\phi_t}(\xi)=e^{-t|\xi|^{\alpha}},\qquad
\widehat{u}(\xi,t)=e^{-t|\xi|^{\alpha}}\widehat{f}(\xi).
\]

By Shannon’s sampling theorem, every \(f\in PW_c\) is stably reconstructible from uniform samples
\(\{f(k/\Omega): k\in\mathbb{Z}\}\) if and only if the sampling rate \(\Omega\) is at least the
\emph{Nyquist rate} \(c/\pi\) (equivalently, the spacing satisfies \(1/\Omega\le \pi/c\)).
However, for uniform spatial sampling, it is not possible to reduce the spatial rate below Nyquist
by augmenting with \emph{all} time samples \(f_t\): in a representative case (e.g., \(\alpha=1\)),
there exists a bandlimited signal $f=f_0$ such that $\|f_0\|=1$ and  $\int_0^1 \sum_{k\in\Z}\left|f_t\left(\frac{m\pi}{c}k\right)\right|^2\,\mathrm{d}t$ is arbitrarily small
 whenever \(\Omega<c/\pi\); hence, a stable reconstruction fails from such uniform space-time data $\{\frac{m\pi}{c}k\}\times [0,1]$.

A remedy is to adopt \emph{periodic non-uniform} spatial patterns \(\Lambda\subset\mathbb{R}\),
which do yield a meaningful \emph{space-time trade-off} \cite{ADK15}: one can choose \(\Lambda\) of
sub-Nyquist density so that the semi-continuous frame inequality holds,
\begin{equation}\label{SemiContFr}
\mathfrak m \,\|f\|_2^2 \;\le\; \int_{0}^{1}\ \sum_{\lambda\in\Lambda} \bigl|f_t(\lambda)\bigr|^2\,dt
\;\le\; \mathfrak M\,\|f\|_2^2,\qquad \forall\, f\in PW_c,
\end{equation}
for some constants \(0<\mathfrak m\le \mathfrak M<\infty\).

The achievable trade-off is limited by the desired numerical conditioning. For instance, the next
result links the \emph{maximal gap} of a stable sampling set \(\Lambda\) to the frame bounds in
\eqref{SemiContFr}.

\begin{theorem}\cite{AGHJKR21} \label{thm_gap_intro}
Assume \(\Lambda\subset\mathbb{R}\) satisfies \eqref{SemiContFr}. Then there exists an absolute constant
\(K>0\) such that, for every \(a\in\mathbb{R}\) and every
\[
R \;\ge\; K\,\max\!\left(\frac{\mathfrak M}{\mathfrak m},\,\frac{1}{c}\right),
\]
the interval \([a-R,a+R]\) contains at least one point of \(\Lambda\).
Equivalently, every gap of \(\Lambda\) has length at most \(2R\).
\end{theorem}
In terms of \emph{Beurling densities}, this yields the bounds
\[
D^{-}(\Lambda)\ \ge\ K^{-1}\,\min\!\left(\frac{\mathfrak m}{\mathfrak M},\,c\right),
\qquad
D^{+}(\Lambda)\ \le\ K\,\mathfrak M,
\]
quantifying how spatial density must scale with the conditioning parameters \({\mathfrak m}, {\mathfrak M}\) and the
band-limit \(c\). 
A more general version of the above result, providing  a more
explicit dependence of $K$ on the parameters of the problem can be found in \cite[Theorem 4]{AGHJKR21}.

Beyond the constraints implied by Theorem~\ref{thm_gap_intro}, the special sampling configurations of Lu and Vetterli that yield \eqref{SemiContFr}, as well as those in \cite{ADK15}, lack the simplicity of regular sampling patterns.
However, it is possible to consider sub-Nyquist equispaced spatial sampling patterns \eqref{SemiContFr} with $\triangle x=\displaystyle \frac{c}{m\pi}$,
$m \in \mathbb{N}$, and restrict the sampling/reconstruction problem to a subset $V \subseteq PW_c$, aiming for an inequality of the form:
\begin{equation}
\label{eq_diff_samp}
\mathfrak m\|f\|_2^2 \le \int_0^1 \sum_{k\in\Z}\left|f_t\left(\frac{m\pi}{c}k\right)\right|^2\,\mathrm{d}t \le \mathfrak M\|f\|_2^2,
\qquad f \in V.
\end{equation}
Specifically, we consider the following signal models.

One can identify a set $\E$ with measure arbitrarily close to $1$ such that \eqref {eq_diff_samp} holds with $V=V_\E=\{f\in PW_c:  \supp  \widehat f \subseteq \E\}$. In effect, $\E$ is the set 
$[-c,c]\setminus\mathcal{O}$ where $\mathcal{O}$ is a small open neighborhood of a  finite set, i.e., $\E$ avoids a certain number of ``blind spots.'' 
 \begin{theorem}[\cite{AGHJKR21}] \label{sufconT0}
Let $ \phi \in \Phi_c$ and $m\geq 2$ be an integer. Then for any  $r>0$ there exists a  compact set $\E\subseteq [-c,c]$ of measure at least $2c-r$ such that %
 any $f\in V_\E$ can be recovered from the samples 
$$
\mathcal M =\left\{f_t\left(\frac{m \pi}{c}k\right): k \in \mathbb{Z}, 0 \leq t \leq 1\right\}
$$
in a stable way. %
\end{theorem}

The set $\E$ in the above theorem depends only on $\phi$ and the choice of $r$. The stable recovery in this case means that \eqref {eq_diff_samp} holds with $\mathfrak M=1$ and some $\mathfrak m > 0$ which is estimated in
a more explicit version of the above result,  
 \cite[Theorem 2.8]{AGHJKR21}. Some results on space-time sampling along curvilinear lattices and their connection to frames can be found in \cite{Zlo22}.

Additional work on discrete- and continuous-time dynamical frames includes
operator-theoretic and structural characterizations of frames generated by
orbits, such as commuting and multi-operator constructions, spectral
conditions, and model-space representations
\cite{BHKLL25,ACCP21,ACCP23,CRS26,CH22_2,Men22}.

Further contributions address multivariate, tensor-product, and
higher-dimensional extensions of dynamical frame constructions
\cite{ZL19,ZLL17_1,BC25,SR25,Cor21,NV23,SCM25}.

Related operator-theoretical techniques for various classes of linear operators in different classes of Hilbert and Banach spaces are studied in \cite{AJ25,BSS25,AJS25,BSK24}.

Finally, we mention frames by iterations in band-limited functions and shift-invariant spaces, which appear in \cite{AADP13,HX22,Zlo22,UZ21}.

\section{Source recovery and frames} \label {SR}
\subsection{Discrete-time model for source recovery}
\label{DTSR}

In this setting, the dynamics are given by \eqref{DS}. The goal is to recover the unknown forcing term
\(w(n)\in\mathcal{W}\subset\HH\) for each \(n\in\mathbb{N}\).
The case \(\mathcal{W}\subsetneq\HH\) is considerably richer and more challenging than the ambient case
\(\mathcal{W}=\HH\).
A motivating application is environmental monitoring, where one seeks to identify the locations and
time-varying magnitudes of pollutant sources from sparse measurements.
We deploy a finite sensor set \(\G=\{g_j\}_{j\in J}\) and collect space-time data \(y_{n,j}\) at sensor \(j\) and time \(n\).

A concrete example is air-quality monitoring as in \cite{RDCV12}, where the objective is to estimate emissions from a
known number of smokestacks using a limited number of sensors.
Here the source subspace can be modeled as
\[
\mathcal{W}=\operatorname{span}\{\psi_1,\dots,\psi_r\},\qquad
w(n)=\sum_{i=1}^r a_i(n)\,\psi_i,
\]
with \(\psi_i\) encoding fixed source locations and \(a_i(n)\) their unknown amplitudes.
Recovering \(w(n)\) from the measurements \(\{y_{n,j}\}\) yields the emission profiles \(\{a_i(n)\}\).
In many practical scenarios, a periodic structure is reasonable, e.g., diurnal behavior:
\(a_i(n+N)=a_i(n)\) for some period \(N\).

It is convenient to arrange the data as a matrix \(Y=[y_{n,j}]\).
Over one period, we view such matrices as elements of a Banach space \(\mathcal{M}_N\subset\mathbb{C}^{N\times J}\) endowed with a suitable norm.
Likewise, an \(N\)-periodic source \(w\) is naturally identified with a vector in \(\HH^N=\underbrace{\HH\times\cdots\times\HH}_{N\ \text{times}}\) with the product inner product.

\medskip
\noindent\textbf{Problem statement.}
Let \(A:\HH\to\HH\) be bounded with \(\|A\|<1\), and let \(\mathcal{W}\subseteq\HH\) be a closed subspace.
Consider the discrete dynamical system \eqref{DS} with an unknown initial state \(x(0)=x_0\) and an unknown \(N\)-periodic forcing \(w(n+N)=w(n)\) taking values in \(\mathcal{W}\).

\medskip
\noindent Find conditions on \(\G=\{g_j\}_{j\in J}\subset\HH\) under which \(w(n)\) can be recovered from the data
\(\{y_{n,j}\}_{n\ge0,\ j\in J}\) \emph{stably}; that is, there exists a bounded linear operator
\[
R:\mathcal{M}_N\longrightarrow \mathcal{W}^{N},\qquad
R\bigl(\{y_{n,j}\}\bigr)=\bigl(w(0),\dots,w(N-1)\bigr).
\]

\medskip

If \(\HH\) is finite dimensional and \(w(n)\equiv w\) is constant in time, the condition \(\|A\|<1\) is not needed; necessary and sufficient conditions on \(\G=\{g_j\}_{j\in J}\) are given in \cite{ADM24}.
In infinite dimensions, stability is nontrivial; necessary and sufficient conditions on \(\G\) were first established in \cite{ADGMM23} for the case \(w(n)\equiv w\in\mathcal{W}\) (constant in time), within a setting more general than \eqref{DS} that includes nonlinearities as well as higher-order difference equations. For this case, it was also shown that when \(\mathcal{W}=\HH\), the problem becomes structurally simple: a necessary and sufficient condition for recovery is that \(\G=\{g_j\}_{j\in J}\) is a frame for \(\HH\) \cite{ADGMM23}. In this case, one needs only finitely many time samples to reconstruct the source term. By contrast, when \(\mathcal{W}\subsetneq\HH\), the situation is subtler: time sampling may need to be infinite \cite{ADGMM23}. For the proper-subspace case \(\mathcal{W}\subsetneq\HH\), the characterization can be stated in terms of projections of the sensing set \(\G=\{g_j\}_{j\in J}\) onto certain invariant subspaces \cite{ACM25, ADGMM23}.
In the periodic case \(w(n)\in\mathcal{W}\), for the linear dynamical system \eqref{DS}, the situation is described by the following theorem.
\begin{theorem}[\cite {ACM25}]\label{MT}
Let \(A:\HH\to\HH\) be bounded with \(\|A\|<1\), let \(N\in\mathbb{N}\), and let \(\mathcal{W}\) be a closed subspace of \(\HH\).
For fixed \(A\), $N$, and \(\mathcal{W}\),
the following are equivalent:
\begin{enumerate}
\item There exists a bounded operator \(R:\mathcal{M}_N\to\mathcal{W}^N\) such that for each solution \(x\) of \eqref{DS} we have
\(
R\bigl(\{\langle x(n),g_j\rangle\}_{n\ge0,\ j\in J}\bigr)=w.
\)
\item For each \(s=0,\dots,N-1\), the family \(\{P_{\mathcal{W}}(T_s^{*}g_j)\}_{j\in J}\) is a frame for \(\mathcal{W}\), where \(P_{\mathcal{W}}\) denotes the orthogonal projection onto \(\mathcal{W}\), and where $T_s := \left( e^{2 \pi i s/N}I-A\right)^{-1}$.
\end{enumerate}
\end{theorem}
A randomized space-time sampling for recovering source terms can be found in \cite{GH25}.

\subsection {Continuous-time model}\label {CTSR}
In the continuous-time model \eqref{DFM}, the source term is \(f+\eta\), where \(f\) is the signal of interest and
\(\eta\) is a Lipschitz-continuous background term.
The signal \(f\) can be modeled as a finite superposition of causal bursts,
\begin{equation}\label{Force}
  f(t)=\sum_{j=1}^N h_j\,\phi_j(t-t_j)\,\chi_{[t_j,\infty)}(t),
\end{equation}
with \(0<t_1<\cdots<t_N\), amplitudes \(h_j\in\HH\), and kernels \(\phi_j\ge 0\) on \([0,\infty)\) satisfying a prescribed decay.
Thus, recovering \(f\) amounts to estimating the onset times \(\{t_j\}\) and the burst shapes/amplitudes \(\{h_j\}\).

Measurements are acquired at sample times \(t_i\ge 0\) by sensors \(g_j\in G=\{g_j\}_{j\in J}\):
\begin{equation}\label{measurements}
  \mathfrak m(t_i,g_j)\;=\;\bigl\langle u(t_i),\,g_j\bigr\rangle+\nu(t_i,g_j),
\end{equation}
where \(\nu\) models measurement noise and \(u\) solves \eqref{DFM} with source \(f+\eta\). Alternatively, weighted average measurements over small time intervals can be used.
 
Rigorous recovery results and stability guarantees have been established under progressively richer assumptions on the pulse shapes \(\phi_j\):
(i) \emph{impulsive sources} \(\phi_j=\delta\) (Dirac bursts) \cite{AHKK23};
(ii) \emph{exponentially decaying} kernels \cite{AGK23};
(iii) \emph{multi-exponential} bursts with unknown onsets and decay rates,
\(f(t)=\sum_{j=1}^{M} h_j\,e^{-\rho_j(t-t_j)}\,\chi_{[t_j,\infty)}(t)\) \cite{AGKMT25},
where \(t_j,\rho_j\), and \(h_j\in\HH\) are all unknown.

A fundamental idea behind the recovery of the source term $f$ is the structural design of the spatial sensor frame set \(\G=\{g_j\}_{j\in J}\). The sensors are placed at locations \(\{x_j\}_{j\in J}\) and acquire samples on the
space-time grid $\{(x_j,n\beta): j\in J,\ n\in\mathbb{N}\}$, where \(\beta>0\) is the sampling time interval.
With an appropriate design of \(\G\), the space-time data allow us to \emph{detect} whether a burst onset \(t_k\) occurs in
the interval \([n\beta,(n+1)\beta)\).
If an onset is detected in \([n\beta,(n+1)\beta)\), then the samples collected up to time \((n+1)\beta\) can be used to
estimate the remaining unknowns -- either \(\rho_k\) or \(h_k\), depending on the forcing model \(f\).
This yields a reconstruction algorithm with guarantees that depend on the sampling interval \(\beta\), the noise level
\(\sigma\), the Lipschitz constant \(L\) of the background term \(\eta\), and the relevant model parameters.

The simplest instance arises for \emph{impulsive} sources (\(\phi_j=\delta\)), where a concrete choice of \(\G\) together
with a simple predictor for \(t_k\) and a direct reconstruction of \(h_k\) provide a complete recovery of
\begin {equation} \label {DM}
f(t)=\sum_j h_j\,\delta(t-t_j).
\end{equation}
For example, the following statement summarizes one of the main theorems in  \cite{AHKK23}.
\begin{theorem}[\cite{AHKK23}] \label{cor:stab_rec_dirct} 
A frame set $\G$ for 
$\HH$ can be designed such that, under appropriate assumptions, a function \(\tilde f(t)=\sum_j \tilde h_j\,\delta(t-\tilde t_j)\) can be constructed from the samples \eqref {measurements} that well approximate the function $f$ in \eqref {DM}.
 In particular, $|t_j-\tilde t_j| \le \frac \beta 2$ 
\begin{equation}
    \label{Herr}
{\|\tilde h_j - h_j\|} \le C_1 L R\beta^2+  C_2\sigma+D(\beta)R\|h_j\|,
\end{equation}
where, $\sigma=\sup\limits_{t,g} |\nu(t,g)|$, $L$ the Lipschitz constant for the background $\eta$, and $R=\sup\limits_{g \in  \G} \|g\|.$
\end{theorem}

\section {System Identification}\label {SI}
System identification is arguably the most intriguing and, at the same time, the least developed category of problems within dynamical sampling. Here, one seeks to identify or estimate various features of the dynamical system \eqref{DS} encoded by the parameters of the evolution operator $E$. In more mathematical terms, the problem is to recover the spectrum of the operator $E$ from dynamical samples $y_{n,j}$. In this section, we will primarily describe the finite-dimensional results obtained in the original paper on the subject \cite{AK16} and their infinite-dimensional version presented in \cite{Tan17}. Related numerical results can also be found in \cite{AHKLLV18}. Finally, we will briefly mention methodologically related results from \cite{KP25, LMT21}. Additionally, an interested reader may want to check out \cite{AGKMT25, BH23, ACCP23, AP23, CT22, ZLL17}.

The investigation of system identification in dynamical sampling started with the following special case of model \eqref{DS}:
\begin{equation}\label{DSinit}
    y_\ell = S_\Omega E^\ell x, \quad x\in \CC^d, \ \ell\in\N_0,
\end{equation}
where $S_\Omega$ is the subsampling operator as in \eqref{e:sds}.

In the case when $\Omega = \Z_d$, i.e.~$S_\Omega = I$, there are classical methods for determining the spectrum $\sigma(E)$ from the values $x$, $Ex$, $E^2x$, \ldots, which are based on the analysis of Krylov subspaces
\[
\mathcal{K}_r(E,x) = {\mathrm{span} }\, \{E^k x: k = 0, 1, \ldots, r \}.
\]
In \cite{AK16}, such a method was generalized to include the case $\Omega\subsetneq \Z_d$; in fact, it can be used for a general linear operator $A$ in place of $S_\Omega$. The generalization is based on the analysis of the \emph{altered} Krylov subspaces
\[
\mathcal{AK}_{r,A}(E,x) = A(\mathcal{K}_r(E,x)) =  {\mathrm{span} }\, \{AE^k x: k = 0, 1, \ldots, r \}.
\]
In particular, one needs to solve the linear equations
\begin{equation}
    \label{gensysid}
    AE^{r+k} x + \sum^{r-1}_{\ell=0} \alpha_\ell AE^{\ell+k} x = 0, \quad k = 0, \ldots, r-1,
\end{equation}
for the minimal value of $r$ for which a solution $\alpha = (\alpha_0, \ldots, \alpha_{r-1})$ exists. The spectrum $\sigma(E)$ will then necessarily contain the roots of the polynomial $$p_\alpha(t) = t^r + \alpha_{r-1}t^{r-1}+\ldots+ \alpha_0.$$
As a consequence, the following result was obtained.

\begin{theorem}[\cite{AK16}, Theorem 3.6]\label{thm:prony}
        Let $E$ be an unknown operator and the set $\Omega \subseteq \{0, 1, \ldots, d-1\}$ be such that any $x\in \CC^d$ could be recovered from the dynamical samples in \eqref{DSinit} with $\ell = 0, \ldots, s-1$, for some $s\in\N$ if the operator $E$ were known. Then $\sigma(E)$ can be recovered from the same dynamical samples with $\ell = 0, \ldots, 2s-1$, for almost every $x \in \CC^d$.
\end{theorem}

If the set $\Omega$ in the above theorem is too small to allow full recovery of the spectrum $\sigma(E)$, then a subset of the so-called observable eigenvalues would be recovered from the dynamical samples. 

In the special case when $\Omega = \{0,m, 2m, \ldots, d-m\}$ for some $m$ that divides $d$ and $E$ is a cyclic convolution operator with a filter $e \in \mathbb C^d$: 
\begin{equation}
    \label{e:convo}
    (Ex)(n) = (e*x)(n) = \sum_{k \in\Z_d} e(k)x(n-k), \ n\in\Z_d,
\end{equation}
the system \eqref{gensysid} reduces to $\frac dm$ systems
   \begin{equation}
       \label{E:ps}
\sum\limits_{\ell = 0}^{r_j-1} \widehat y_{k+l}(j)\alpha_\ell(j) = -\widehat y_{k+r_j} (j) , \quad
k = 0, \ldots, m - 1, \quad j = 1, \ldots, d/m.
   \end{equation}
Moreover, if $m=d$ and $E$ is the circular shift, $(Ex)(n)=(Ex)(n+1)$, $n\in\Z_d$, then the solution of \eqref{E:ps} yields the coefficients of the classical Prony polynomial \cite{P795, FR13}, which is used to solve the famous problem of identifying an $s$-sparse vector from $2s$ of its consecutive Fourier coefficients. In fact, in the latter case, a straightforward application of Theorem \ref{thm:prony} yields the original Prony algorithm.

In \cite{Tan17}, Tang followed the same approach for the case of convolution operators in $\ell^2(\Z)$, i.e.~$Ex = e*x$ with $e, x \in\ell^2(\Z)$. In fact, equations \eqref{E:ps} remain essentially unchanged. The key difference is that $\widehat y$ is then a periodic function rather than a finite-dimensional vector, and one has to solve uncountably many systems similar to \eqref{E:ps}. In practice, it is often assumed that $\supp e$ and $\supp x$ are finite sets, i.e.~both $e$ and $x$ have a finite impulse response. In that case, it suffices to solve finitely many systems to recover both $e$ and $x$.

Numerically, spectrum recovery can be performed using various algorithms.
The original Prony algorithm is typically very sensitive to noise. However, various modifications have been devised to improve stability in the classical case. In \cite{Tan17}, it is shown how the SVD-based matrix pencil method and the ESPRIT method generalize to the setting of dynamical sampling with convolution operators. 

\begin{remark}
    Observe that the matrices of the linear systems used for spectral identification, such as \eqref{E:ps}, are Hankel. This allows one to employ Cadzow denoising techniques when implementing 
    spectral recovery from noisy data \cite{AHKLLV18, Tan17}, see also \cite{CHLY25}. Numerical simulations described in \cite{AHKLLV18, Tan17} show that even the straightforward recovery algorithm becomes fairly robust with respect to additive measurement noise when coupled with Cadzow denoising. 
\end{remark}

The approach to system identification developed in \cite{AK16, Tan17} has been expanded in various directions.
For example, it allowed one to view Prony's problem as the spectral identification problem for an unknown restriction of a known linear operator (a similar point of view was adopted in \cite{PP13}).  
In \cite{KP25}, the authors developed this idea in the context of Banach modules. Using \eqref{gensysid}, they obtained, for example, an algorithm for identifying finite linear combinations of time-frequency shift operators. In \cite{LMT21}, related ideas were used for identification of homogeneous polynomials to enable learning of interaction kernels in various dynamical systems.

\section{Outlook and open problems}\label{sec:outlook} \label {outlook}

Dynamical sampling has matured into a broad interface between harmonic analysis, operator theory, approximation theory, complex analysis, and control theory. Even so, a number of important questions remain open -- ranging from foundational existence and structure problems for frames of iterations, to quantitative sampling design, to algorithmic guarantees in source recovery and system identification. We highlight several promising directions for future research.

\subsection*{(A) Frames of iterations: existence, structure, and indices}
\begin{itemize}
  \item \textbf{Normal vs.~self-adjoint operators.} For a bounded self-adjoint operator $A$ in an infinite dimensional Hilbert space, the normalized orbits $\{A^n g/\|A^n g\|\}_{g \in \G,\;n\geq 0}$ never form a frame \cite{ACMT17}. The normal case remains open. We state the following
conjecture:
\begin{conjecture}
If $A$ is a bounded normal operator on an infinite dimensional Hilbert space $\HH$ and  $\mathcal{G} \subset \HH,$ then the system  $\left\{\frac{A^n g}{\norm{A^n g}}\right\}_{g \in \G,\;n\geq 0}$ is not a frame for $\HH$.
\end{conjecture}
 
  \item \textbf{Redundancy of Carleson frames.} Carleson frames exhibit striking ``excessive redundancy.'' The full extent of this redundancy, however, is unknown. Conjecture \ref{conj:carl} illustrates our expectations.
  \item \textbf{Frame index and minimal generators.} A formula for the Parseval index $\gamma_p(T)$ of a bounded operator $T$ 
on a Hilbert space $\HH$ has been established as 
$
\gamma_p(T) = \dim \overline{(I - TT^*)(\HH)}.
$
The general frame index $\gamma(T)$ can be obtained by minimizing the 
Parseval index over all contractions similar to $T$ whose adjoint is 
strongly stable. An intriguing question is whether one can find 
\emph{intrinsic}, purely operator-theoretic expressions for $\gamma(T)$.
Furthermore, how does this index behave when additional structure is 
imposed -- such as when $T$ is subnormal, or pure contraction of a particular class?
  \item \textbf{From orbits to multi-orbits.} Beyond single-operator orbits, the multi-operator setting  is only partly charted (commuting or nearly commuting families have been considered). A comprehensive structure theory for frames generated by several operators, with useful quantitative bounds on the number of generators and iterates, is open.
\end{itemize}

\subsection*{(B) Semi-continuous frames and time discretization}
\begin{itemize}
  \item \textbf{Discretization on unbounded windows.} Theorems~\ref{ScToDscr}  bridge continuous and discrete-time on $[0,L]$. What persists on $[0,\infty)$?
  \emph{If a semi-continuous frame holds on $[0,\infty)$, does there exist a \underline{discrete-time} set with finite upper Beurling density that yields a frame?} (Cf.\ Open Problem~\ref{OP:L-infty}.)
  \item \textbf{Normal reductive case.} The stability‐for‐all‐$L$ phenomenon (Theorem~\ref{SCFrSA}) is proved for invertible self-adjoint $A$; the normal reductive case remains a conjecture 
  (see Conjecture \ref{conj:reduct}).
  \item \textbf{Mobile dynamical sampling.} Creating a theory for semi-continuous frames of the form $\{A^tg_j(t)\}_{t\in[0,L], j\in J}$ would lead to better understanding of mobile sampling in the context of dynamical systems.
\end{itemize}

\subsection*{(C) Developing space-time trade-off}
\begin{itemize}
  \item \textbf{Nonuniform spatial grids.} Obtaining quantitative estimates for space-time trade-off on nonuniform temporal or spatial grids would lead to better practical sampling designs.     
  A key investigative frontier is to \emph{optimize} the sampling patterns in terms of their stability, i.e.~frame bounds.
  \item \textbf{Model-restricted recovery.} Sub-Nyquist equispaced sampling becomes viable on structured subclasses $V\subset PW_c$ (``blind-spot'' avoidance \cite{AGHJKR21}). Characterizing  the \emph{largest} admissible subclasses for a given step $\Delta x$, specific kernels $\varphi_t$, and fixed frame bounds remains open.
\end{itemize}

\subsection*{(D) Source recovery: robustness and design}
\begin{itemize}
  \item \textbf{Stability bounds.} Existing stability bounds in source recovery problems \cite{AHKK23, AGK23, AGKMT25} are conservative in the sense that they pertain to the worst-case scenario. Estimating average stability bounds for a given sampling design and noise level is a wide open problem, for which even a proper formulation has not been settled. 
  \item \textbf{Optimal sampling scheme design.} 
  In \cite{AHKK23, AGK23, AGKMT25}, reconstruction guarantees of predictive algorithms depend on various parameters such as the sampling schedule and the structure and location of the sensors. Designing \emph{near‐optimal} sensor configurations and sampling schedules that minimize error constants while respecting budgets is a ripe optimization problem with immediate impact.
  
  \item\textbf{Burst detection with background.} Predictive algorithms of \cite{AHKK23, AGK23, AGKMT25} are effective in separating the source of interest from background because these two types of external sources are very different from each other (e.g.~burst-like sources vs.~Lipschitz background). What other types of sources could be separated through similar algorithms?
  
  \item \textbf{Stability and design in discrete-time cases.} The stability of the reconstruction operator \(R\) in Theorem~\ref{MT} or related Theorems in \cite{ADGMM23} depends on the sensing set \(\G=\{g_j\}\).
Quantifying this stability dependence and
designing \(\G\) optimally under a fixed cardinality constraint are important problems.
  \end{itemize}

\subsection*{(E) System identification from dynamical samples}
\begin{itemize}
  \item \textbf{Practical algorithms for system identification.} The theory for system identification developed in \cite{AK16, Tan17, KP25} calls for the development of better numerical methods in various scenarios, as well as adapting existing numerical methods to wider contexts. For example, can matrix-pencil and ESPRIT-type methods employed in \cite{Tan17} be used in the framework of \cite{KP25}?
  \item \textbf{Nonlinear measurements.} System identification from non-linear measurements (e.g.\ phaseless \cite{BH23}, saturated \cite{AFGJR25}, or fractional part \cite{BKP22}) is a wide open problem. A similar question can be posed for other core dynamical sampling problems, and the relevant research is only just beginning to emerge \cite{AKT18, BH23, Boz24, BK23, ADGMM23}. 
\end{itemize} 

\subsection*{(F) Dynamical sampling on graphs and in Banach spaces}
\begin{itemize}
  \item \textbf{Dynamical sampling on graphs.} The study of core dynamical sampling problems in a graph-theoretical setting has only just begun \cite{ABKMP24, HNT24, GH25, HLTY25}. There are many exciting questions that one can pursue in transferring the results from static signal processing on graphs to the dynamical setting.
  \item \textbf{Banach space setting.} Many dynamical sampling problems, including system identification and construction of dynamical frames, may naturally be posed in the setting of Banach spaces \cite{KP25}. However, the systematic extension of the dynamical sampling theory from Hilbert to Banach spaces remains largely unexplored. 
\end{itemize}

Progress on these questions would clarify the geometry of dynamical frames, sharpen our understanding of the interplay between the continuous-time and discrete-time settings of dynamical sampling, and strengthen our ability to design and analyze algorithms for source recovery and system identification.

\section*{Acknowledgment}
The authors thank the reviewers for their helpful comments and suggestions.

\vspace{1cm}
On behalf of all authors, the corresponding author states that there is no conflict of interest.

\bibliographystyle{siam}

\end{document}